\def\sqr#1#2{{\vcenter{\vbox{\hrule  height.#2pt
	\hbox{\vrule width.#2pt height#1pt \kern#1pt \vrule width.#2pt}
	\hrule height.#2pt}}}}
\def\sq{\sqr55}        %    A small square for end-of-proofs. 
\documentclass[11pt]{article}
\usepackage{amssymb,amsmath,theorem,rlepsf}
\newcommand{\stdspace}{\hskip 0.75em plus 0.15em\ignorespaces}
\newenvironment{proof}{\par{\bf Proof}\stdspace}{\hfill$\sq$\par}
\renewenvironment{abstract}{\par{\bf Abstract}\par}{\par}
\newcommand{\qed}{\hfill$\sq$\par}
\newcommand{\co}{\colon\thinspace}    %  Colon with correct spacing for maps.
\newcommand{\ppar}{\par\vskip 8pt plus4pt minus4pt}

\newcommand{\rk}[1]{\ppar{\bf #1}\stdspace}
\let\fonote\footnote
\def\footnote#1{\fonote{\small #1}}

\def\thanks#1{\relax}
\def\theoremstyle#1{\relax}
\def\dedicatory#1{\relax}
\def\tableofcontents{\relax}
\def\maketitle{\relax}
\def\title#1{\def\thetitle{#1}}
\def\authors#1{\def\theauthors{#1}}
\def\author#1{\def\theauthors{#1}}
\def\address#1{\def\theaddress{#1}}

\def\email#1{\def\theemail{#1}}
\def\url#1{\def\theurl{#1}}
\long\def\abstract#1\endabstract{\long\def\theabstract{#1}}
\def\primaryclass#1{\def\theprimaryclass{#1}}
\def\secondaryclass#1{\def\thesecondaryclass{#1}}
\def\keywords#1{\def\thekeywords{#1}}
\theorembodyfont{\sl}  %%% default gt theorembody font
\makeatletter
\long\def\@caption#1[#2]#3{\par\addcontentsline{\csname
  ext@#1\endcsname}{#1}{\protect\numberline{\csname
  the#1\endcsname}{\ignorespaces #2}}\begingroup
    \@parboxrestore
    \small
    \@makecaption{\csname fnum@#1\endcsname}{\ignorespaces #3}\par
  \endgroup}
\makeatother
\newtheorem{theorem}{Theorem}[section]
\newtheorem{lemma}[theorem]{Lemma}
\newtheorem{proposition}[theorem]{Proposition}

\newtheorem{claim}{Claim}
\theoremstyle{definition}
\theorembodyfont{\rm}
\newtheorem{definition}[theorem]{Definition}

\begin{document}
%                          Title page  
%
\input gtoutput
\volumenumber{1}\papernumber{3}\volumeyear{1997}
\pagenumbers{21}{40}\published{30 July 1997}

\title{Canonical Decompositions of 3--Manifolds}
\shorttitle{Canonical decompositions of 3--manifolds}

\authors{Walter D Neumann\\Gadde A Swarup}

\address{Department of Mathematics\\The University of Melbourne\\Parkville, Vic 3052\\Australia\\}

\email{neumann@ms.unimelb.edu.au{\rm\ \  and\ \ }gadde@ms.unimelb.edu.au}
\let\theurl\relax

\abstract
We describe a new approach to the canonical decompositions 
of 3--manifolds along tori and annuli due to
Jaco--Shalen and Johannson (with ideas from Waldhausen) ---
the so-called JSJ--decomposition theorem.
This approach gives an accessible proof of the decomposition
theorem;
in particular it does not use the annulus--torus
theorems, and the theory of Seifert fibrations does not need to be
developed in advance.
\endabstract

\primaryclass{57N10, 57M99}

\secondaryclass{57M35}

\keywords{3--manifold, torus decomposition, JSJ--decomposition, 
Seifert\break \hbox spread 3pt{manifold, simple manifold}}

\proposed{David Gabai}\received{25 February 1997}
\seconded{Robion Kirby, Ronald Stern}\accepted{27 July 1997}
\maketitlepage

%%%%%%%%%%%%%%%%%%%%   End of title page
%%%%%%%%%%%%%%%%%%%%    Start of main body of article

\section{Introduction}

\thanks{This research is supported by the Australian Research Council}

In this paper we describe a proof of the so-called JSJ--decomposition
theorem for 3--manifolds.  This proof was developed as an exercise for
ourselves, to confirm an approach that we hope will be useful for 
JSJ--decomposition in the group-theoretic context. It seems to give a
particularly accessible proof of JSJ for 3--manifolds, involving few
prerequisites.  For example, it does not use the annulus--torus
theorems, and the theory of Seifert fibrations does not need to be
developed in advance.

We do not give the simplest version of our proof.  We prove JSJ for
orientable Haken 3--manifolds with incompressible boundary. This
involves splitting along tori and annuli.  As we describe later, if
one restricts to the case that all boundary components are tori, then
one only needs to split along tori.  With this restriction our proof
becomes very much simpler, with many fewer case distinctions.  The
general JSJ--decomposition can then be deduced quite easily from this
special case by doubling the 3--manifold along its boundary.  We took
the more direct but less simple approach because this was an
exercise to test a concept: there is no clear analogue of the
peripheral structure and that of a double in the group theoretic case
and we preferred our approach to be closer to possible group theoretic
analogues.

We do not prove all the properties of the JSJ--decomposition, although
it seems that this could be done from our approach with a little more
effort. The main results we do not prove --- the enclosing theorem and
Johannson deformation theorem --- have very readable accounts in
Jaco's book \cite{jaco}.

\rk{Acknowledgement}This research was supported by the Australian 
Research Council.

\section{Canonical Surfaces and W--Decomposition}

We consider orientable Haken 3--manifolds with incompressible boundary
and consider splittings along essential annuli and tori.  Let $S$ be
an annulus or torus that is properly embedded in $(M,\partial M)$ and
which is essential (incompressible and not boundary-parallel).

\begin{definition}
$S$ will be called \emph{canonical} if any other properly 
embedded essential annulus or torus $T$ can be isotoped
to be disjoint from $S$.

Take a disjoint collection $\{S_{1},\dots,S_{s}\}$ of canonical 
surfaces in $M$ such that 
\begin{itemize}
\item no two of the $S_{i}$ are parallel;
\item the collection is maximal among disjoint collections of 
canonical surfaces with no two parallel.
\end{itemize}
A maximal system exists because of the Kneser--Haken finiteness
theorem.  The result of splitting $M$ along such a system will be
called a
\emph{Waldhausen decomposition} (or briefly \emph{W--decomposition}) of
$M$.  It is fairly close to the usual JSJ--decomposition which will be
described later. The maximal system of pairwise non-parallel canonical
surfaces will be called a \emph{W--system}.
\end{definition}
The following lemma shows that the W--system $\{S_{1},\dots,S_{s}\}$ is
unique up to isotopy.
\begin{lemma}\label{all at once}
Let $S_1,\ldots,S_k$ be pairwise disjoint and non-parallel canonical
surfaces in $(M,\partial M)$. Then any incompressible annulus or torus
$T$ in $(M,\partial M)$ can be isotoped to be disjoint from
$S_1\cup\dots\cup S_k$. Moreover, if $T$ is not parallel to any $S_i$
then the final position of $T$ in $M-(S_1\cup\dots\cup S_k)$ is
determined up to isotopy.
\end{lemma}
\begin{proof}
By assumption we can isotop $T$ off each $S_i$ individually. Writing
$T=S_0$, the lemma is thus a special case of the following stronger
result:

\begin{lemma}
Suppose $\{S_0,S_1,\ldots,S_k\}$ are incompressible surfaces in
an irreducible manifold $M$ such that each pair can be isotoped to be
disjoint.  Then they can be isotoped to be pairwise disjoint and the
resulting embedded surface $S_0\cup\ldots\cup S_k$ in $M$ is
determined up to isotopy.
\end{lemma}

Without the uniqueness statement a quick conceptual proof of this uses
minimal surfaces (one can also use a traditional cut-and-paste type
proof; see below). Choose a metric on $M$ which is product metric near
$\partial M$. Freedman, Hass and Scott show in
\cite{freedman-hass-scott} that the least area representatives of the
$S_i$ are pairwise disjoint or identical.  When two least area
representatives are identical, we isotop them off each other in a
normal direction (no two of the $S_i$ embed to identical one-sided
least area surfaces since an embedded one-sided surface can never be
isotoped off itself: the transverse intersection with an isotopic copy
is a $1$--manifold dual to the first Stiefel--Whitney class of the
normal bundle of the surface).  There is a slight complication in that
the least area representative $S_i\to S_i'\subset M$ may not be an
embedding; it can be a double cover onto an embedded one-sided surface
$S_i'$. But in this case we use a nearby embedding onto the boundary
of a tubular neighbourhood of $S_i'$.

In \cite{freedman-hass-scott} surfaces are considered up to homotopy
rather than isotopy, so the above argument uses the fact that two homotopic
embeddings of an incompressible surface are isotopic (see
\cite{waldhausen annals}, Corollary 5.5).

For the uniqueness statement assume we have $S_1,\ldots,S_k$
disjointly embedded and then have two different embeddings of $S=S_0$
disjoint from $T=S_1\cup \ldots\cup S_k$.  Let $f\co S\times I\to
M$ be a homotopy between these two embeddings and make it transverse
to $T$.  The inverse image of $T$ is either empty or a system of
closed surfaces in the interior of $S \times I$. Now use Dehn's Lemma
and Loop Theorem to make these incompressible and, of course, at the
same time modify the homotopy (this procedure is described in Lemma
1.1 of \cite{waldhausen} for example).  We eliminate 2--spheres in
the inverse image of $T$ similarly.  If we end up with nothing in the
inverse image of $T$ we are done.  Otherwise each component $T'$ in
the inverse image is a parallel copy of $S$ in $S\times I$ whose
fundamental group maps injectively into that of some component $S_i$
of $T$. This implies that $S$ can be homotoped into $S_i$ and its
fundamental group $\pi_1(S)$ is conjugate into some $\pi_1(S_i)$.  It
is a standard fact (see eg \cite{swarup indian})
in this situation of two incompressible surfaces having
comparable fundamental groups that, up to conjugation, either
$\pi_1(S)=\pi_1(S_j)$ or $S_j$ is one-sided and $\pi_1(S)$ is the
fundamental group of the boundary of a regular neighbourhood of $T$
and thus of index $2$ in $\pi_1(S_j)$. We thus see that either $S$ is
parallel to $S_j$ and is being isotoped across $S_j$ or it is a
neighbourhood boundary of a one-sided $S_j$ and is being isotoped
across $S_j$.  The uniqueness statement thus follows.

One can take a similar approach to prove the existence of the isotopy
using Waldhausen's classification \cite{waldhausen annals} of proper
incompressible surfaces in $S\times I$ to show that $S_0$ can be
isotoped off all of $S_1,\dots,S_k$ if it can be isotoped off each of
them.
\end{proof}
%\begin{ssect}
The thing that makes decomposition along incompressible annuli and
tori special is the fact that they have particularly simple
intersection with other incompressible surfaces.
\begin{lemma}
If a properly embedded incompressible torus or annulus $T$ in an
irreducible manifold $M$ has been isotoped to intersect another
properly embedded incompressible surface $F$ with as few components in
the intersection as possible, then the intersection consists of a
family of parallel essential simple closed curves on $T$ or possibly a
family of parallel transverse intervals if $T$ is an annulus.
\end{lemma}
\begin{proof}
We just prove the torus case.  Suppose the intersection is
non-empty. If we cut $T$ along the intersection curves then the
conclusion to be proved is that $T$ is cut into annuli.  Since the
Euler characteristics of the pieces of $T$ must add to the Euler
characteristic of $T$, which is zero, if not all the pieces are annuli
then there must be at least one disk.  The boundary curve of this disk
bounds a disk in $F$ by incompressibility of $F$, and these two disks
bound a ball in $M$ by irreducibility of $M$.  We can isotop over this
ball to reduce the number of intersection components, contradicting
minimality.
\end{proof}
%\end{ssect}
\section{Properties of the W--decomposition}
%\begin{ssect}
Let $M_1,\ldots,M_m$ be the result of performing the W--decomposition
of $M$ along the W--system $\{S_1\cup\dots\cup S_s\}$.

We denote by $\partial_1M_i$ the part of $\partial M_i$ coming from
$S_1\cup\dots\cup S_s$ and by $\partial_0M_i$ the part coming from
$\partial M$.  Thus $\partial_0M_i$ and $\partial_1M_i$ are
complementary in $\partial M_i$ except for meeting along circles.  The
components of $\partial_1M_i$ are annuli and tori. Both
$\partial_1M_i$ and $\partial_0M_i$ are incompressible since we
started with an $M$ with $\partial M$ incompressible.  However, it is
possible that $\partial M_i$ is not incompressible.

\begin{figure}[htbp]
\centerline{
\relabelbox\small
\epsfxsize.8\hsize\epsffile{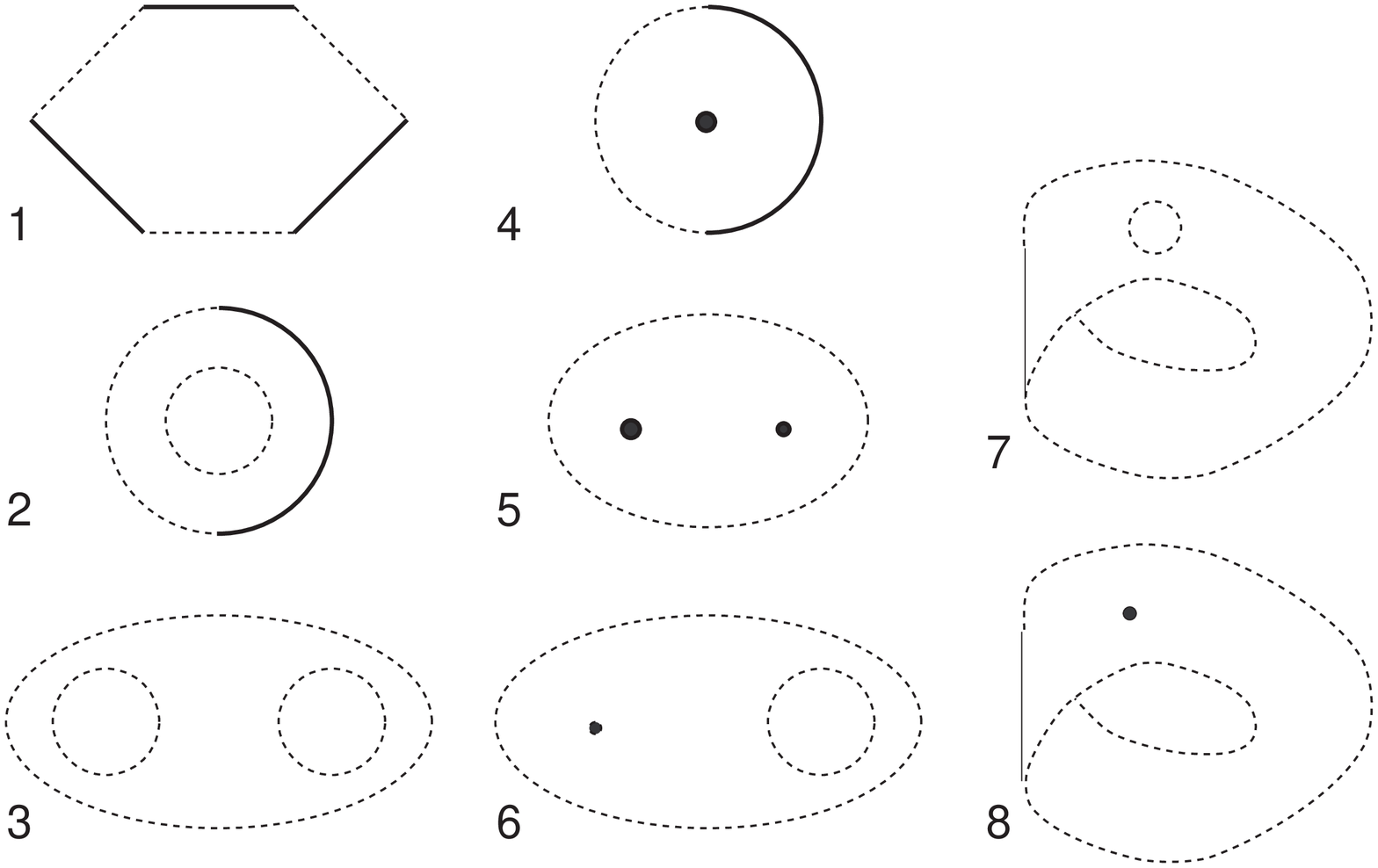}
\relabel{1}{1}
\relabel{2}{2}
\relabel{3}{3}
\relabel{4}{4}
\relabel{5}{5}
\relabel{6}{6}
\relabel{7}{7}
\relabel{8}{8}
\endrelabelbox
}
\caption{\small \label{eight}The eight Seifert fibred building blocks.  
In each case we have drawn the base surface with the portion of the
boundary in $\partial_0$ resp.\ $\partial_1$ drawn solid resp.\
dashed. The dots in cases 4,5,6 represent singular fibres and in
case 8 it represents a fibre that may or may not be singular. These
cases thus represent infinite families.}
\end{figure}

Lemma \ref{all at once} shows that any incompressible annulus or torus
$S$ in $M$ can be isotoped into the interior of one of the $M_i$.  If
this $S$ is an annulus its boundary $\partial S$ is necessarily in
$\partial_0M_i$. By the maximality of $\{S_1,\dots,S_s\}$, if $S$ is
canonical it must be parallel to one of the $S_j$.  However, there may
be such $S$ which are not canonical and not parallel to an $S_j$.
%\end{ssect}
\begin{definition}
We call $M_i$ \emph{simple} if any essential annulus or torus
$(S,\partial S)\subset(M_i,\partial_0M_i)$ is parallel to
$\partial_1M_i$.  We call $M_i$ \emph{special simple} if, in addition,
it admits an essential annulus in $(M_i,\partial_1 M_i)$ (allowing the
possibility that the annulus may be parallel to $\partial_0M_i$).
This will, in particular, be true if some component of $\partial_0M_i$
is an annulus.
\end{definition}
\begin{proposition}\label{nonsimple is seifert}
If $M_i$ is non-simple then $(M_i,\partial_0M_i)$ is either Seifert
fibred or an $I$--bundle.  

If $M_i$ is special simple then either $(M_i,\partial_0M_i)$ is
Seifert fibred of one of the types illustrated in Figure\ \ref{eight} or
$(M_i,\partial_0 M_i)=(T^2\times I,\emptyset)$ (in which case $M$ is a
$T^2$--bundle over $S^1$ with holonomy of trace
$\ne\pm2$).
\end{proposition} 
A consequence of this proposition is that an irreducible manifold $M$
that has an essential annulus or torus but no canonical one is an
$I$--bundle or Seifert fibered.
\begin{proof}
We first consider the non-simple case.  There are a few cases that
will be exceptional from the point of view of our argument and that we
will need to treat specially as they come up.  They are the cases 
\begin{itemize}
\item $M$ is an $I$--bundle over a torus or Klein bottle --- then $M$
also fibres as an $S^1$--bundle over annulus or M\"obius band
respectively;
\item $M$ is an $S^1$--bundle over torus or Klein bottle.
\end{itemize}
These are, in fact, among the non-simple manifolds that admit more
than one fibred structure. In these cases it is not hard to see that
no essential annulus or torus is canonical, so the W--decomposition of
$M$ is trivial, so $M=M_1$.

We drop the index and denote
$(N,\partial_0N,\partial_1N)= (M_i,\partial_0M_i,\partial_1M_i)$.

Consider a maximal disjoint collection of pairwise non-parallel
essential annuli and tori $\{T_1,\dots,T_r\}$ in
$(N,\partial_0N)$. Split $N$ along this collection into pieces
$N_1,\dots,N_n$.  We still denote by $\partial_0N_i$ the part coming
from $\partial M$, that is, $\partial_0N_i=\partial_0N\cap
N_i=\partial M\cap N_i$, and we denote by $\partial_1N_i$ the rest of
$\partial N_i$.

We shall show that if $M$ is not one of the exceptions listed above,
then the pieces $N_i$ are fibred either as $I$--bundles or Seifert
fibrations, and, moreover, these fibrations match up when we glue the
$N_i$ together to form $N$.  In fact:
\begin{claim}\label{claim1}
Each $N_i$ is either an
$I$--bundle over a twice punctured disk, a M\"obius band, or a
punctured M\"obius band or is Seifert fibred of one of the types in
Figure\ \ref{eight}.
\end{claim}

$\partial_1N_i$ consists of annuli and tori, some of which may come
from the $S_j$, but at least one of which comes from a $T_j$.  For
this $T_j$ we let $T'_j$ be an essential annulus or torus that
intersects $T_j$ essentially and we assume the intersection $T_j\cap
T'_j$ is minimal.

The
proof is a case by case analysis of this situation.  Since the proofs
in the different cases are fairly similar, we give a complete argument
in only a couple of typical cases.

{\bf Case A}\stdspace $\partial T_j\cap\partial T'_j\ne\emptyset$.  We shall
see that $(N_i,\partial_0N_i)$ is an $(I,\partial I)$--bundle over a
surface and $\partial_1N_i$ is the part of the bundle lying over the
boundary of the surface.  The surface in question is either a twice
punctured disk, a M\"obius band, or a punctured M\"obius band (except
when $N=M$ was itself an $I$--bundle over a torus or Klein bottle).

This will be proved case by case.  Denote a component of $T'_j\cap
N_i$ which intersects $T_j$ by $P$ and let $s$ be a component of the
intersection $P\cap T_j$.  Since $T'_j\cap(T_1\cup\dots\cup T_t)$ will
be a finite number of segments crossing $T'_j$, there will be another
segment $s'$ in $\partial P\cap(T_1\cup\dots\cup T_t)$ and the rest of
$\partial P$ will consist of two segments in $\partial_0N_i$.  Let
$s'$ be in $T_k$.  We distinguish two subcases.

{\bf Case A1}\stdspace $T_j\ne T_k$. Note that $P$ is an $I$--bundle over an
interval with $s$ and $s'$ the fibres over the ends of this interval,
and $T_j$ and $T_k$ are $I$--bundles over circles.  Cutting $T_j$ and
$T_k$ along $s$ and $s'$, pushing them just inside $N_i$, and then
pasting parallel copies of $P$ gives an annulus (see Figure\
\ref{figA1}). 
\begin{figure}[htbp]
\centerline{\relabelbox\small
\epsfxsize.5\hsize\epsffile{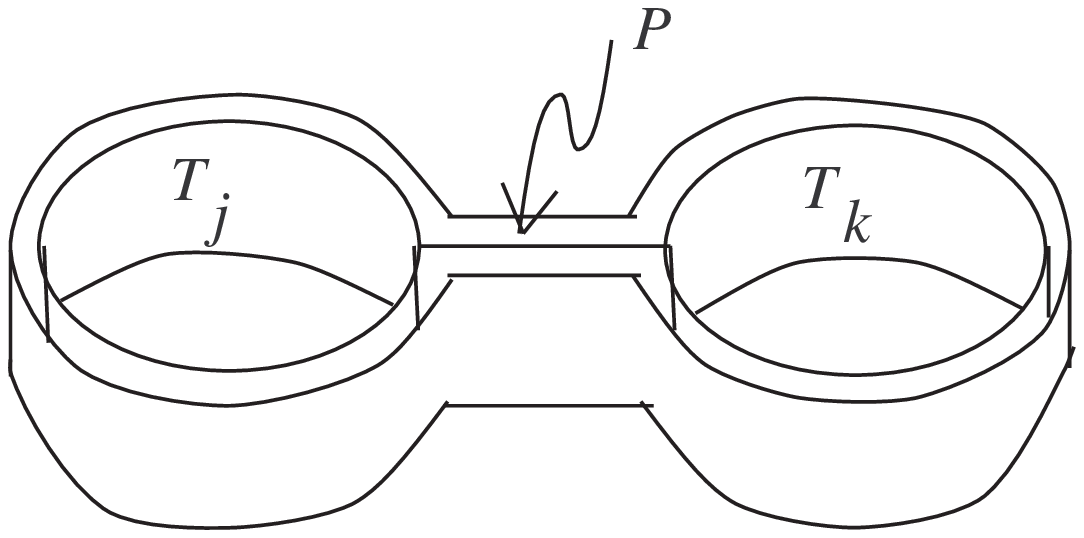}
\adjustrelabel <-6pt, 5pt> {P}{$P$}
\adjustrelabel <-2pt, 5pt> {j}{$T_j$}
\adjustrelabel <-2pt, 5pt> {k}{$T_k$}
\endrelabelbox}
\caption{\small \label{figA1}Case A1}
\end{figure}
It cannot be parallel into $\partial_0N_i$ because then $T_j$ would be
inessential (note that our assumption of boundary-incompressibility
implies that an annulus is inessential if an arc across the annulus is
isotopic into the boundary). It can not be isotopically trivial since
then $T_j$ and $T_k$ would be parallel.  It is thus parallel into
$\partial_1N_i$.  Thus $N_i$ is an $I$--bundle over a twice punctured
disk.

{\bf Case A2}\stdspace $T_j=T_k$. Since $P$ is an $I$--bundle over an interval
with $s$ and $s'$ the fibres over the ends of this interval, and $T_j$
is an $I$--bundle over a circle, we may orient the fibres of these
bundles so that $s$ is oriented the same in each.  Then $s'$ may or
may not be oriented the same in each. Moreover, $P$ may meet $T_j$ at
$s'$ at the same side of $T_j$ as it meets it at $s$ or at the
opposite side.

{\bf Case A2.1.1}\stdspace $P$ meets $T_j$ both times from the same side
and $s'$ has the same orientation in both $I$--bundles.
Then if we cut $T_j$ along $s\cup s'$ and glue two parallel copies
$P'$ and $P''$ of $P$ we get two annuli in $N_i$ (see Figure\
\ref{figA211}).
\begin{figure}[htbp]
\centerline{\relabelbox\small
\epsfxsize.5\hsize\epsffile{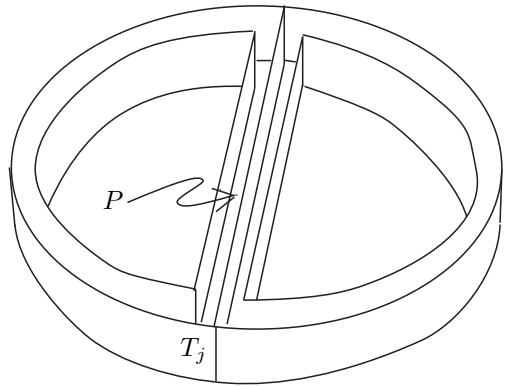}
\adjustrelabel <-2pt, 6pt> {P}{$P$}
\adjustrelabel <-6pt, 4pt> {j}{$T_j$}
\endrelabelbox}
\caption{\small \label{figA211}Case A2.1.1}
\end{figure}
If either of these annuli is trivial in $N_i$ we could
have removed the intersection $s\cup s'$ of $T_j$ and $T'_j$ by an
isotopy.  Thus they are both non-trivial and must be parallel to
components of $\partial_1N_i$ or $\partial_0N_i$.  If either is
parallel into $\partial_0N_i$ then $T_j$ would have been
inessential.  They are thus both parallel into $\partial_1N_i$ and
we see that $N_i$ is an $I$--bundle over a twice punctured disk.

{\bf Case A2.1.2}\stdspace $P$ meets $T_j$ both times from the same side and
$s'$ has opposite orientations in the two $I$--bundles. Then cutting
$T_j$ as above and gluing two copies $P'$ and $P''$ of $P$ yields a
single annulus (see Figure\ \ref{figA212}) 
\begin{figure}[htbp]
\centerline{%
\relabelbox\small
\epsfxsize.5\hsize\epsffile{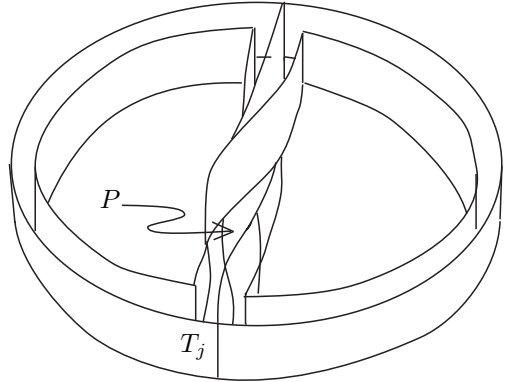}
\adjustrelabel <-3pt, 5pt> {P}{$P$}
\adjustrelabel <-6pt, 4pt> {j}{$T_j$}
\endrelabelbox
}
\caption{\small \label{figA212}Case A2.1.2}
\end{figure}
which may be isotopically
trivial or parallel into $\partial_1N_i$ (as before, it cannot be
parallel into $\partial_0N_i$).  This gives an $I$--bundle over a
M\"obius band or punctured M\"obius band.

{\bf Case A2.2}\stdspace $P$ meets $T_j$ from opposite sides (and $s'$ has the
same or opposite orientations in the two $I$--bundles). After cutting
open along $T_j$ we have two copies $T_j^{(1)}$ and $T_j^{(2)}$ of
$T_j$ and this becomes similar to case A1: cutting $T_j^{(1)}$ and
$T_j^{(2)}$ along $s$ and $s'$ and pasting parallel copies of $P$
gives an annulus. As in Case A1, this annulus cannot be parallel into
$\partial_0N_i$. It may be parallel into $\partial_1N_i$, and then
$N_i$ is an $I$--bundle over a twice punctured disk.  Unlike case A1,
it may also be isotopically trivial, in which case $N_i$ is an
$I$--bundle over an annulus, so $N=M$ is an $I$--bundle over the torus
or Klein bottle, giving two of the exceptional cases mentioned at the
start of the proof.

{\bf Case B}\stdspace $\partial T_j\cap \partial T'_j=\emptyset$. In this case
we will see that $N_i$ is Seifert fibred of one of eight basic types
of Figure \ref{eight} (in the exceptional cases mentioned at the start
of the proof --- where the W--decomposition is trivial and $N=M$ is
itself an $S^1$--bundle over an annulus, M\"obius band, torus, or Klein
bottle --- there is just one piece $N_1$ after cutting and it does not
occur among the eight types).

Denote a component of $T'_j\cap N_i$ which intersects $T_j$ by $P$ and
let $s$ be a component of the intersection $P\cap T_j$. This $P$ is an
annulus and $s$ is one of its boundary components.  Denote the other
boundary component by $s'$.  It either lies on $\partial_0N_i$ or on
some $T_k$.

{\bf Case B1}\stdspace $s'$ lies on $\partial_0N_i$. There are two subcases
according as $T_j$ is an annulus or torus. 

{\bf Case B1.1}\stdspace $T_j$ is an annulus (and $s'$ lies on
$\partial_0N_i$; here, and in the following, subcases inherit the
assumptions of their parent case).  Cutting $T_j$ along $s$ and
pasting parallel copies of $P$ to the resulting annuli gives a pair of
annuli which must both be parallel into $\partial_1N_i$, since if
either were parallel into $\partial_0N_i$ we could have removed the
intersection $s$ of $T_j$ and $T'_j$ by an isotopy.  This gives Figure\
\ref{eight}.1 (ie case 1 of Figure\ \ref{eight}).

{\bf Case B1.2}\stdspace $T_j$ is a torus.  Cutting $T_j$ along $s$ and
pasting parallel copies of $P$ to the resulting annulus gives an
annulus which must be parallel into $\partial_1N_i$, since if it were
parallel into $\partial_0N_i$ then $T_j$ would be boundary parallel.
This gives Figure\ \ref{eight}.2.

{\bf Case B2}\stdspace $s'$ lies on $T_k$ and $T_k\ne T_j$. There are three
subcases according as both, one, or neither of $T_j$ and $T_k$ are
annuli.

{\bf Case B2.1}\stdspace $T_j$ and $T_k$ both annuli.  Cutting $T_j$
and $T_k$ along $s$ and $s'$ and pasting parallel copies of $P$ to the
resulting annuli gives a pair of annuli.  They cannot both be parallel
into $\partial_1N_i$, for then we would get a picture like Figure\
\ref{eight}.1 but based on an octagon rather than a hexagon and an
annulus joining opposite $\partial_0$ sides of this would contradict
the fact that the $T_i$'s formed a maximal family. They cannot both be
parallel into $\partial_0N_i$ for then $T_i$ and $T_j$ would be
parallel.  Thus one is parallel into $\partial_1N_i$ and the other
into $\partial_0N_i$.  This gives Figure\ \ref{eight}.1. again.

{\bf Case B2.2}\stdspace One of $T_j$ and $T_k$ an annulus and the other a
torus.  Cut and paste as before gives a single annulus which may be
parallel into $\partial_0$, giving Figure\ \ref{eight}.2.  It cannot be
parallel into $\partial_1$, since then we would get a picture like
Figure\ \ref{eight}.2 but with two $\partial_0$ annuli in the outer
boundary rather than just one. One could span an essential annulus
from one of these components of $\partial_0$ to itself around the
torus, contradicting maximality of the family of $T_j$'s.

{\bf Case B2.3}\stdspace Both of $T_j$ and $T_k$ tori.  Cut and paste as
before gives a single torus.  It cannot be parallel into $\partial_0$
since then the annulus that can be spanned across $P$ from this torus
to itself would contradict the fact that the $T_j$'s form a maximal
family.  It may be parallel into $\partial_1$ leading to Figure\
\ref{eight}.3 or it may also bound
a solid torus which gives \ref{eight}.6.

{\bf Case B3}\stdspace $s'$ lies on $T_j$. There two
subcases according as $T_j$ is an annulus or torus, and each subcase
splits into three subcases according to whether $P$ meets $T_j$ from
the same or opposite sides at $s$ and $s'$, and if the same side, then
whether fibre orientations match at $s$ and $s'$. We shall describe
them briefly and leave details to the reader.

{\bf Case B3.1}\stdspace $T_j$ an annulus.

{\bf Case B3.1.1}\stdspace $P$ meets the annulus $T_j$ both times
from the same side.

{\bf Case B3.1.1.1}\stdspace The orientations of $s$ and $s'$ match.  Cut and
paste as before gives an annulus and a torus.  The annulus cannot be
parallel into $\partial_1$ since otherwise one can find an annulus
that contradicts maximality of the family of $T_j$'s, so it is parallel
into $\partial_0$.  Similarly, the torus cannot be parallel into
$\partial_0$. This leads to the cases of Figure\ \ref{eight}.2 or
\ref{eight}.4 according to whether the torus is parallel into
$\partial_1$ or bounds a solid torus.

{\bf Case B3.1.1.2}\stdspace Orientations of $s$ and $s'$ do not match.  Cut
and paste as before gives an annulus.  Whether it is parallel into
$\partial_0$ or $\partial_1$, $N_i$ is a circle bundle over a M\"obius
band with part of its boundary in $\partial_0$ and an annulus from
$\partial_0$ to $\partial_0$ running once around the M\"obius band
contradicts the maximality of the family of $T_j$'s.  Thus this case
cannot occur. 

{\bf Case B3.1.2}\stdspace Assumptions of Case B3.1 and $P$ meets the annulus
$T_j$ from opposite sides. This gives the same as Case B2.1, except
for one extra possibility, since both annuli of Case B2.1 being
parallel into $\partial_0$ is not ruled out.  This extra case leads to
$N=M$ being an $S^1$--bundle over the annulus or M\"obius band.

{\bf Case B3.2}\stdspace Assumptions of Case B3 with $T_j$ a torus.

{\bf Case B3.2.1}\stdspace $P$ meets the torus $T_j$ both times
from the same side.  

{\bf Case B3.2.1.1}\stdspace Orientations of $s$ and $s'$
match.  Cut and paste gives two tori, neither of which
is parallel to $\partial_0$.  This leads to Figure\
\ref{eight}.3, \ref{eight}.5,
 or \ref{eight}.6.  

{\bf Case B3.2.1.2}\stdspace Orientations of $s$ and $s'$ do not
match.  Cut and paste gives one torus which may be parallel to
$\partial_1$ or bound a solid torus. This gives Figure\
\ref{eight}.7 or \ref{eight}.8.

{\bf Case B3.2.2}\stdspace Assumptions of Case B3.2 and $P$ meets the torus
$T_j$ from opposite sides. This gives the same as Case B2.3, except
that in the situation of Figure\ \ref{eight}.6 the ``singular'' fibre
need not actually be singular, in which case $N=M$ is an $S^1$--bundle
over the torus or Klein bottle (these are among the special manifolds
listed at the start of the proof).

This completes the analysis of the pieces $N_i$ and thus proves 
Claim \ref{claim1}. 
We must now verify that the fibred structures on the pieces
$N_i$ match up in $N$.

Suppose two pieces $N_{i_1}$ and $N_{i_2}$ meet across $T_j$ and let
$T'_j$ be as above.  The above argument showed that a fibred
structure can be chosen on these two pieces to match the fibred
structure of the $T'_j$ and hence to match each other across $T_j$.
Thus if every piece $N_i$ has a unique fibred structure then the
fibred structures must match across every $T_j$.

\begin{claim}\label{claim2}
 The only pieces $N_i$ for which the fibration is not
unique up to isotopy are:

{\bf1\rm)}~~$I$--bundle over M\"obius band, which also admits a fibration as

{\bf2\rm)}~~the Seifert fibration of Figure\ \ref{eight}.4 with degree
$2$ exceptional fibre;

{\bf3\rm)}~~the Seifert fibration of Figure\ \ref{eight}.5 with both
exceptional fibres of degree 2, which also admits the structure of

{\bf4\rm)}~~the circle bundle of Figure\ \ref{eight}.8 with no exceptional
fibre.
\end{claim}

If we grant this claim then the matching of fibrations on the various
$N_i$ follows: since in the above non-unique cases $N_i$ has only one
boundary component $T_j$, we simply choose the fibration on this $N_i$
that is forced by $T'_j$ and it will then match the neighbouring
piece.

Claim \ref{claim2} follows by showing that the fibration is unique in all
other cases.  It is clear that the only $N_i$ that is both an
$I$--bundle and Seifert fibred is the one of cases 1 and 2, so we need
only consider non-uniqueness of Seifert fibration.  All we really need
for the above proof is that the Seifert fibration is unique up to
isotopy when restricted to the boundary, which is clear for cases 1,
2, 4 of Figure\ \ref{eight} and follows from the fact that a circle
fibre generates an infinite cyclic normal subgroup of the fundamental
group in the other cases. The uniqueness of the fibration on the whole
of $N_i$ follows, if desired, by a standard argument once at least one
fibre has been determined.

To complete the proof of the proposition we must consider the case
that $M_i$ is special simple, so it is simple but has an essential
annulus in $(M_i,\partial_1M_i)$.  If we call this annulus $P$ we can
repeat word for word the arguments of Case B above to see that $M_i$
is of one of the eight types listed in Figure\ \ref{eight} or is an
$S^1$--bundle over an annulus or M\"obius band with boundary belonging
to $\partial_1M_i$. Since $S^1$--bundle over M\"obius band is included
in the eight types and $S^1$--bundle over annulus is $T^2\times I$,
the statement of the proposition follows.  Note that the $T^2\times I$
case can only occur if the two boundary tori are the same in $M$, in
which case $M$ is a torus bundle over the circle.  The holonomy of
this torus bundle cannot have trace $\pm2$ since then the torus would
not be canonical (in fact $M$ would be an $S^1$--bundle over torus or
Klein bottle).
\end{proof}
%\begin{ssect}
It is worth noting that the above proof shows also that the only
non-simple Seifert fibred manifolds with non-unique Seifert fibration
are those with trivial W--decomposition mentioned at the beginning of
the above proof (circle bundles over annulus, M\"obius band, torus, or
Klein bottle) and the one that comes from matching the non-unique
cases 3 and 4 of the above Claim \ref{claim2} together, giving the Seifert
fibration over the projective plane with unnormalized Seifert
invariant $\{-1;(2,1),(2,-1)\}$ (two exceptional fibres of degree 2
and rational Euler number of the fibration equal to zero).  This
manifold has two distinct fibrations of this type.  Note that matching
two of case 3 or two of case 4 together gives the tangent circle
bundle of the Klein bottle, which is one of the examples with trivial
W--decomposition already mentioned, but we see this way its Seifert
fibration over $S^2$ with four degree two exceptional fibres
(unnormalized invariant $\{0;(2,1),(2,1),(2,-1),(2,-1)\}$).
%\end{ssect}

%\begin{ssect}
We shall classify the pieces $M_i$
with $\partial_1\ne\emptyset$ into three types:
\begin{itemize}
\item $M_i$ is \emph{strongly simple}, by which we mean simple and 
not special simple;
\item $M_i$ is an $I$--bundle;
\item $M_i$ is Seifert fibred.
\end{itemize}
We first discuss which $M_i$ belong to more than one type.
%\end{ssect}
\begin{proposition}
The only cases of an $M_i$ having more than one type are:

$\bullet$~~The $I$--bundle over M\"obius band is also Seifert fibered.

$\bullet$~~The $I$--bundles over twice punctured disk and once
punctured M\"obius band are also strongly simple.
\end{proposition}
\begin{proof}
If $M_i$ admits both a Seifert fibration and an $I$--bundle structure,
it is an $I$--bundle over an annulus or M\"obius band (since we are
assuming $\partial_1\ne\emptyset$). If $M_i$ is $I$--bundle over
annulus then $\partial_1M_i$ consists of two annuli.  Since they are
parallel in $M_i$ they must have been equal in $M$, so $M$ is obtained
by pasting these two annuli together, that is, $M$ is an $I$--bundle
over torus or Klein bottle.  But the annulus is then not canonical (in
fact, such $M$ has trivial W--decomposition), so this cannot occur.
Thus the only case is that $M$ is $I$--bundle over M\"obius band.

If $M_i$ is both strongly simple and an $I$--bundle then we have
already pointed out that $M_i$ cannot be $I$--bundle over annulus.  The
$I$--bundle over M\"obius band is Seifert fibred and is therefore
special simple (see below).  An $I$--bundle over a bounded surface of
Euler number less than $-1$ will not be simple. Thus the only cases
are those of the proposition.

Suppose $M_i$ is Seifert fibred and simple. If some component of
$\partial_0M_i$ is an annulus then $M_i$ is special simple. Otherwise,
$\partial_1M_i$ consists of tori and $M_i$ is Seifert fibred with at
least two singular fibres if the base is a disk and at least one if it
is an annulus.  It is thus easy to see that for each component of
$\partial_1M_i$ there will be an essential annulus with both
boundaries in this component, so $M_i$ is again special simple.  Thus
$M_i$ cannot be both Seifert fibred and strongly simple.
\end{proof}
%\begin{ssect}
We now want to investigate the possibility of adjacent fibred pieces
in the W--decomposition having matching fibres where they
meet.  We therefore assume that the W--decomposition is non-trivial, so
$M_i$ has non-empty $\partial_1$.  

Suppose the pieces $M_i$ and $M_k$ on the two sides of a canonical
annulus $S_j$ are both $I$--bundles.  If $i\ne k$ then in each of $M_i$
and $M_k$ we can find an essential $I\times I$ from $S_j$ to itself.
Gluing these gives an essential annulus crossing $S_j$ essentially.
Thus $S_j$ was not canonical.  If $i=k$ the argument is similar. Thus
$I$--bundles cannot be adjacent in the W--decomposition.
%\end{ssect}

%\begin{ssect}
It remains to consider the case that both pieces on each side of a
canonical annulus or torus $S_j$ are Seifert fibred.  Suppose the
fibrations on both sides match along $S_j$. 

If a piece $M_i$ adjacent to $S_j$ is not simple we shall decompose it
as in the proof of Proposition \ref{nonsimple is seifert} and, for the
moment, just consider the simple piece of $M_i$ adjacent to $S_j$.  We
will call it $N_i$ for convenience.  We first assume that $N_i$ is not
pasted to itself across $S_j$.

This $N_i$ has an embedded essential annulus compatible with its
fibration of one of the following types:
\begin{itemize}
\item
an essential annulus from $S_j$ to a component of $\partial_0N_i$
(``essential'' here means incompressible and not parallel into
$\partial_1N_i$ by an isotopy that keeps the one boundary component in
$\partial_0$ and the other in $\partial_1$, but of course allows the
first to isotop to $\partial_0\cap\partial_1$);
\item 
an essential annulus from $S_j$ to $S_j$ in $N_i$
(incompressible and not parallel to $\partial_1N_i$).
\end{itemize}
This can be seen on a case by case basis by considering the eight
cases of Figure\ \ref{eight}.

If we had an essential annulus of the first type in $N_i$ then,
depending on the situation on the other side of $S_j$, we can either
glue it to a similar annulus on the other side to obtain an essential
annulus in $(M,\partial M)$ crossing $S_j$ in a circle, or glue two
parallel copies to an essential annulus from $S_j$ to itself on the
other side, to obtain an essential annulus crossing $S_j$ in two
circles. Either way we see that $S_j$ was not canonical, so this
cannot occur.

Thus the only possibility is that we have essential annuli from $S_j$
to itself on both sides of $S_j$, which we can then glue together to
get a torus or Klein bottle crossing $S_j$ in two circles.  If it were
a Klein bottle, then the boundary of a regular neighbourhood would be
an essential torus crossing $S_j$ in four circles, contradicting that
$S_j$ is canonical. If it is a torus it will be essential unless it is
parallel into $\partial M$.  The only way this can happen is if the
pieces on each side of $S_j$ are of the type of case 2 or 4 of Figure\
\ref{eight} and $S_j$ is the annular part of $\partial_1$ of these
pieces.  Moreover, we see that the $M_i$ on each side of $S_j$ is
simple, for $N_i$ would otherwise have to be case 2 of Figure\
\ref{eight} with another Seifert building block pasted along the
inside torus component of $\partial_1N_i$ and we could find an annulus
in $M_i$ from $S_j$ to itself that goes through this additional
building block and is therefore not parallel to $\partial M$.

\begin{figure}[htbp]
\centerline{\epsfxsize.8\hsize\epsffile{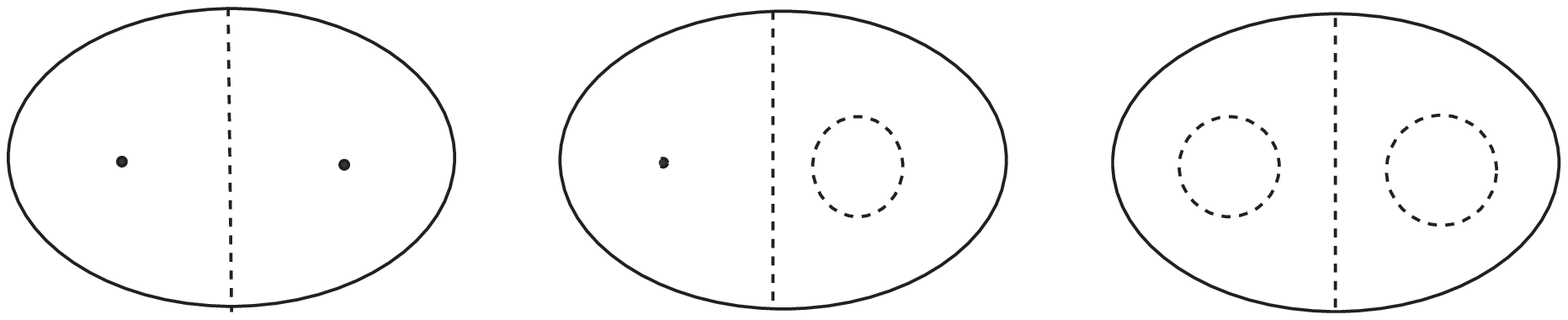}}
\centerline{\epsfxsize.17\hsize\epsffile{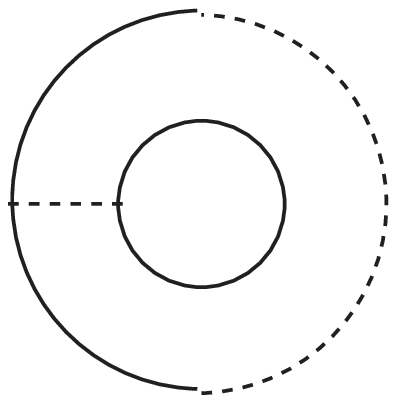}}
\caption{\small \label{extra annuli fig}Matched annuli.  In each case we have
drawn the base surface of the Seifert fibration with the portion of
the boundary in $\partial_0$ drawn solid. The matched annulus lies over
the dashed line in the interior.}
\end{figure}

We must finally consider the case that $N_i$ is pasted to itself
across $S_j$.  An annulus connecting the two boundary components of
$N_i$ that are pasted becomes a torus or Klein bottle in $M$.  If it
were a Klein bottle, then the boundary of a regular neighbourhood
would be an essential torus crossing $S_j$ in two circles,
contradicting that $S_j$ is canonical. If it is a torus it is
essential unless it is parallel into $\partial M$.  By considering the
cases in Figure\ \ref{eight} one sees that the only possible $N_i$ is
Figure\ \ref{eight}.1 pasted as in the bottom picture of Figure\
\ref{extra annuli fig}.  Moreover, as in the previous paragraph we see
that there cannot be another Seifert building block pasted to the
remaining free annulus.

Summarising, we have:
%\end{ssect}
\begin{lemma}\label{extra annuli}
If two fibred pieces of the W--decomposition are adjacent in $M$ with
matching fibrations then they are each of type 2 or 5 of Figure\
\ref{eight} matched along the annular part of $\partial_1$.  If a
fibred piece is adjacent to itself then it is of type 1 of Figure\
\ref{eight}. The possibilities are drawn in Figure\ \ref{extra annuli
fig}.  \qed\end{lemma}
\begin{definition}
We shall call an annulus separating two matched fibred pieces or a
fibred piece from itself as in
Lemma \ref{extra annuli} a \emph{matched annulus}. If we delete all matched
annuli from the W--system then the remaining surfaces will be
called the \emph{JSJ--system} and the decomposition of $M$ along this
JSJ--system will be called the \emph{JSJ--decomposition of $M$}.
\end{definition}

Let $\{S_1,\dots,S_t\}$ be the JSJ--system of canonical surfaces, so
that splitting along this system gives the JSJ--decomposition.  We
recover the W--decomposition as follows: for each piece $M_i$ of the
JSJ--decomposition as in Figure\ \ref{extra annuli fig} we add an annulus
as in that figure.

%\section{Relationship with other versions of \\JSJ--decomposition}
\section{Other versions of JSJ--decomposition}

The JSJ--decomposition is the same as the canonical splitting as described by Jaco
and Shalen in \cite{jaco-shalen} Chapter V, section 4.  That
decomposition is characterised by the fact that it is a decomposition
along a minimal family of essential annuli and tori that decompose $M$
into fibred and simple non-fibred pieces.  The JSJ--system of surfaces
satisfies this minimality condition by construction.

Jaco and Shalen's ``characteristic submanifold'' is a fibred
submanifold $\Sigma$ of $M$ which is essentially the union of the
fibred parts of the JSJ--splitting except that:
\begin{itemize}
\item
wherever two fibred parts of the JSJ--decomposition meet along
an essential torus or annulus, thicken that torus or annulus and add
the resulting $T^2\times I$ or $A\times I$ to the complementary part
$\overline{M-\Sigma}$;
\item wherever two non-fibred pieces meet along an annulus or torus,
thicken that surface and add the resulting $A\times I$ or
$T^2\times I$ to $\Sigma$.
\end{itemize}
%\begin{ssect}
The special case that $M$ has Euler characteristic $0$, so its boundary
consists only of tori, deserves mention.  If there is an annulus in
the maximal system of canonical surfaces, then any adjacent piece
$M_i$ of the W--decomposition has an annular component in its
$\partial_0$ and is therefore Seifert fibred by Proposition
\ref{nonsimple is seifert}. Moreover, the fibrations on these adjacent
pieces match along this annulus, so it is a matched annulus.  Hence:
%\end{ssect}
\begin{proposition}\label{10.6.2}
If $M$ has Euler characteristic $0$ (equivalently, $\partial M$
consists only of tori) then the JSJ--system consists only of tori.
\qed
\end{proposition} 
This proposition seems less well-known than it should be.  It is
contained in Proposition 10.6.2 of \cite{johannson} and is used
extensively without reference in \cite{eisenbud-neumann}.

%\begin{ssect}
There is another modification of the JSJ--decomposition that is often
useful, called the \emph{geometric decomposition}, since it is the
decomposition that underlies Thurston's geometrisation conjecture (or
rather ``geometrisation theorem'' since Thurston proved it in the
Haken case).  The geometric decomposition may have incompressible
M\"obius bands and Klein bottles as well as annuli and tori in the
splitting surface.  
It is obtained from the JSJ--decomposition by eliminating all
pieces which are fibred over the M\"obius band as follows:
\begin{itemize}
\item delete the canonical annulus that bounds any piece which is an
$I$--bundle over a M\"obius band and replace it by the $I$--bundle
over the core circle of the M\"obius band (which is itself a M\"obius
band), and
\item delete the canonical torus that bounds any piece which is an
$S^1$--bundle over a M\"obius band and replace it by the $S^1$--bundle
over the core circle of the M\"obius band (which is a Klein bottle).
\end{itemize}

One advantage of the geometric decomposition is that it lifts 
correctly in
finite covering spaces.
%\end{ssect}

\begin{figure}[htbp]
\centerline{\epsfxsize.17\hsize\epsffile{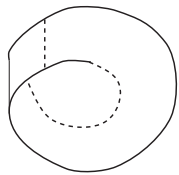}}
\caption{\small \label{annulus-canonical}The annulus over the dashed interval
is annulus-canonical although it has an essential transverse
torus. The dashed portion of the M\"obius band boundary
represents a canonical annulus.}
\end{figure}

\begin{figure}[htbp]
\centerline{\epsfxsize.17\hsize\epsffile{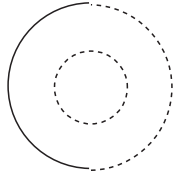}}
\caption{\small \label{torus-canonical}The torus over the inner circle is
always torus-canonical, unless it bounds a solid torus. The dashed
portion of the outer boundary represents a canonical annulus.}
\end{figure}

%\begin{ssect}
We may also consider only essential annuli or only essential tori and
restrict our definition of ``canonical'' to the one type of surface.
Thus we define an essential annulus to be ``annulus-canonical'' if it
can be isotoped off any other essential annulus, and an essential
torus is ``torus-canonical'' if it can be isotoped off any other
essential torus.
%\end{ssect}
\begin{proposition}
An essential annulus is annulus-canonical if and only if it is
canonical or is as in Figure\ \ref{annulus-canonical}.  An essential
torus is torus-canonical if and only if it is canonical or is parallel
to a torus formed from a canonical annulus and an annulus in $\partial
M$ (Figure\ \ref{torus-canonical}).
\end{proposition} 
\begin{proof} The ``if'' is clear.  For the only if, note that 
any annulus-canonical annulus that is not canonical will intersect an
essential torus and must therefore occur inside a Seifert fibred
piece.  It will thus play a similar role to the matched canonical
annuli discussed above.  But the argument we used to classify matched
canonical annuli shows that any matched annulus-canonical annulus is
as in Lemma \ref{extra annuli} and is thus canonical, except that the
argument using the boundary of a regular neighbourhood of a Klein
bottle no longer eliminates the case of Figure\
\ref{annulus-canonical}. The argument in the torus case is similar and
left to the reader.
\end{proof}

\iffalse It is worth noting that if we generalise as in Jaco Shalen 
\cite{jaco-shalen} to consider splittings with respect to a 
submanifold of $\partial M$ rather than all of $\partial M$ then we 
\emph{can} find examples with torus-canonical tori that are not canonical. 
An example is $M=A\times S^{1}$ and $\partial_{0}M=X\times S^{1}$ 
with $X$ a union of  two or more disjoint intervals in $\partial A$.
\smallskip\fi

\section{Only torus boundary components}

If $M$ only has torus boundary components then we can do
W--decomposition from the start only using torus-canonical tori.  This
leads to a much simpler proof. There are many fewer cases to consider
and the issue of ``matched annuli'' disappears --- Seifert fibrations
never match across a canonical torus, so the W--decomposition one gets
by this approach is exactly the JSJ--decomposition.  The general case
as described by Jaco and Shalen can then be deduced from this case by a
doubling argument.  If one is only interested in 3--manifold splittings
this is therefore probably the best approach.

We sketch the argument.  If $M$ is not an $S^1$--bundle over torus or
Klein bottle we cut a non-simple piece $N=M_i$ of this ``toral''
W--decomposition into pieces $N_j$ along a maximal system of
disjoint non-parallel essential tori, analogously to the proof of
Proposition \ref{nonsimple is seifert}.  These pieces turn out to be
of nine basic types, namely what one gets by taking the examples of
Figure\ \ref{eight} that have no annuli in $\partial_0$ (cases
3,5,6,7,8), and then allowing some but not all of the boundary
components to be in $\partial_0$ (so, for example, Case 3 splits into
three types according to whether $0$, $1$, or $2$ of the boundary
components are in $\partial_0$.)  One then shows, as before, that the
fibrations match up to give a Seifert fibration of $N$.

 The classification of pieces $M_i$ that are both simple and Seifert
fibered is well known and can also easily be extracted from our
discussion.  The ones with boundary are precisely cases 3,5,6,7,8 of
Figure\ \ref{eight} but with any number of boundary components in
$\partial_0M=\partial M$.  The closed ones are manifolds with finite
$H_1(M)$ that fiber over $S^2$ with at most three exceptional fibers
or over $\mathbb RP^2$ with at most one exceptional fiber.

\section{Bibliographical Remarks}

The theory of characteristic submanifolds of 3--manifolds for Haken 
manifolds with toral boundaries was first outlined by Waldhausen in 
\cite{waldhausen:conf}; 
see also \cite{waldhausen:recent} for his later account of the topic.
It was worked out in this form by Johannson \cite{johannson}.  The
description of the decomposition in terms of annuli and tori was given
in Jaco--Shalen's memoir \cite{jaco-shalen} (see also Scott's paper
\cite{scott} for this description 
as well as a proof of the Enclosing Theorem).  The idea of canonical
surfaces is suggested by the work of Sela and Thurston; a similar idea
was used independently by Leeb and Scott in \cite{leeb-scott} in the
context of non-positively curved manifolds.  Generalizations using
submanifolds of the boundary are discussed in both Jaco--Shalen and
Johannson; the characteristic submanifold that we described here is
with respect to the whole boundary ($T=\partial M$, in the terminology
of
\cite{jaco}).

\end{document}